\pdfoutput=1

\documentclass[11pt]{amsart}
\usepackage[centertags]{amsmath}
\usepackage{amsfonts}
\usepackage{amssymb}
\usepackage{amsthm}
\usepackage{epsfig}
\usepackage{fullpage }
\usepackage{MnSymbol}
\usepackage{subfigure}
\usepackage{color}

\newtheorem{theorem}{Theorem}
\newtheorem{proposition}{Proposition}
\newtheorem{lemma}{Lemma}
\newtheorem{corollary}{Corollary}
\theoremstyle{definition}
\newtheorem{definition}{Definition}

\theoremstyle{remark}
\newtheorem{remark}{Remark}
\newtheorem{example}{Example}

\begin{document}
\newcommand{\A}{\mathbb{A}}
\newcommand{\B}{\mathbb{B}}
\newcommand{\s}{\mathbb{S}}
\newcommand{\T}{\mathbb{T}}
\newcommand{\G}{\mathbb{G}}

\newcommand{\overlay}{\cupdot}
\newcommand{\ar}[1]{\overrightarrow{#1}}

\newcommand{\dotcup}{\ensuremath{\mathaccent\cdot\cup}}
\newcommand{\du}{\oast}

\title[Bipartite partial duals]{Bipartite partial duals and circuits 
in medial graphs}

\author[S.~Huggett]{Stephen Huggett$^*$}
\author[I.~Moffatt]{Iain Moffatt$^\dagger$}

\begin{abstract}
It is well known that a plane graph is Eulerian if and only if its 
geometric dual is bipartite. We extend this result to partial duals 
of plane graphs.  We then characterize all bipartite partial duals 
of a plane graph in terms of oriented circuits in its medial graph. 

\end{abstract}

\thanks{
${\hspace{-1ex}}^*$School of Computing and Mathematics, University of Plymouth, PL4 8AA, Devon, UK.;  \\
${\hspace{.35cm}}$ \texttt{s.huggett@plymouth.ac.uk}}

\thanks{
${\hspace{-1ex}}^\dagger$
Department of Mathematics and Statistics,  University of South Alabama, Mobile, AL 36688, USA;\\
${\hspace{.35cm}}$ \texttt{imoffatt@jaguar1.usouthal.edu}}

\date{\today}


\maketitle

\section{Introduction and statements of results}

The geometric dual, $G^*$, of an embedded graph $G$ is a  fundamental 
construction in graph theory and appears in many places throughout 
mathematics. Motivated by various new constructions in knot theory, 
S.~Chmutov, in \cite{Ch1}, introduced the concept of the partial dual 
of an embedded graph.  Roughly speaking,  a partial dual  is obtained by forming the geometric dual with respect 
to only a subset of edges  of an embedded graph (a formal definition is given 
Subsection~\ref{ss.pddef}). 
Partial duality appears to be a fundamental operation on embedded 
graphs and, although it has only recently been introduced, it has 
found a number of applications in graph theory, topology, and physics 
(see, for example, \cite{Ch1,EMM,HMV,KRVT09,Mo2,Mo3,Mo4,Mo5,VT,VT10}). 
While geometric duality always preserves  the  surface in which a 
graph is embedded, this is not the case for the more general  partial 
duality. For example, if $G$ is a plane graph, then $G^*$ is also a 
plane graph, but a partial dual $G^A$ of $G$ need not be plane. 
Partially dual embedded graphs can have  different topological and 
graph theoretical properties.

Rather than being concerned with the ways in which a graph and its partial dual can differ, here we are interested in how partial duality both preserves and transforms the structure of an embedded graph. In particular, we determine the extent to which partial duality preserves the following classical connection between Eulerian and bipartite plane graphs. 
\begin{theorem}\label{t.intro1}
Let $G$ be a plane graph, then $G$ is Eulerian if and only if its 
dual, $G^{*}$, is bipartite. \qed
\end{theorem}
This theorem is well-known (see, for example,  Theorem~34.4 of \cite{VLW} or Example~10.2.10 of \cite{BM}). 
It is known to hold more generally for binary matroids (see 
\cite{We69} and also \cite{Ja79}), but it does not hold for non-plane graphs (although the geometric dual of a bipartite graph is always Eulerian). Note that since $(G^*)^*=G$, the words bipartite and Eulerian can be in interchanged in Theorem~\ref{t.intro1}.
Here we give the extension of this classical connection between 
Eulerian and bipartite graphs from geometric duality to partial 
duality. We prove: 
\begin{theorem}\label{t.intro2}
Let $G$ be a plane graph and $A\subseteq E(G)$. Then:
\begin{enumerate}
\item $G^A$ is bipartite if and only if  the components of $G|_{A}$  
and $G^*|_{A^c}$ are Eulerian;
\item $G^A$ is Eulerian if and only if  $G|_{A}$  and $G^*|_{A^c}$ 
are bipartite.
\end{enumerate}
\end{theorem}
This result appears in Section~\ref{s.bepd} as Theorem~\ref{t.pde}. 
Its proof requires much more work than the case for geometric duality 
stated in Theorem~\ref{t.intro1}. To prove the result we introduce a new way of 
obtaining the underlying abstract graph of partial dual $G^A$ 
(see  Theorem~\ref{t.newpds} below). The advantage of this construction is 
that it avoids having to construct the embedded graph $G^A$ itself.

Having established the relation between Eulerian and bipartite partial 
duals, we then turn our attention to the problem of determining which 
subsets of edges in a graph give rise to bipartite and Eulerian 
partial duals. That is, given a plane graph $G$, the problem is to characterize the 
subsets $A\subseteq E(G)$ having the property that $G^A$ is bipartite 
or Eulerian. It turns out that this problem is intimately related to 
oriented circuits in the medial graph $G_m$ of $G$. 
We provide the following complete characterization of edge sets that 
lead to bipartite partial duals:
\begin{theorem}
    Let $G$ be a plane graph. Then the partial dual $G^A$ is 
    bipartite if and only if $A$ is the set of $c$-edges arising from 
    an all-crossing direction of $G_m$. 
\end{theorem}
The terminology for this theorem, together with its proof, appears in 
Section~\ref{s.med}. (Figure~\ref{f.m} offers a quick indication of the terminology). We note that Theorem~\ref{t.intro2} is used in 
an essential way  to prove this characterization.  We also find a 
sufficient condition for a set of edges to give rise to an Eulerian 
partial dual in terms of circuits in medial graphs (see 
Corollary~\ref{c.deu}). Some connections of our results with knot theory are discussed in Remark~\ref{r.knots}.

\section{Embedded graphs and duality}

\subsection{Cellularly embedded graphs and ribbon graphs}

We begin with a brief review of embedded graphs and ribbon graphs. 
Note  that we will be concerned with both cellularly and 
non-cellularly embedded graphs.

An {\em embedded  graph} $G=(V(G),E(G)) \subset \Sigma$ is a graph 
drawn on a surface $\Sigma$  in such a way that edges only intersect 
at their ends. The  arcwise-connected components of $\Sigma 
\backslash G$ are called the {\em regions} of $G$. If each of the 
regions of an embedded graph $G$ is homeomorphic to a disc we say 
that $G$ is a {\em cellularly embedded graph}, and its regions are 
called {\em faces}. A {\em plane graph} is a graph that is cellularly 
embedded in the sphere (rather than the plane).

Two embedded graphs, $G\subset \Sigma$ and $G'\subset \Sigma'$, are 
said to be {\em equal} if there is a homeomorphism from  $\Sigma$ 
to $\Sigma'$ that sends $G$ to $G'$. As is common, we will often 
abuse notation and identify an embedded graph with its equivalence 
class under equality.

We will need to work with cellularly embedded graphs which arise as 
subgraphs of cellularly embedded graphs. Accordingly, we will often 
find it  convenient and natural to describe embedded 
graphs as ribbon graphs.
\begin{definition}
A {\em ribbon graph} $G =\left(  V(G),E(G)  \right)$ is a (possibly 
non-orientable) surface with boundary represented as the union of 
two  sets of topological discs: a set $V (G)$ of {\em vertices}, and 
a set of {\em edges} $E (G)$ such that 
\begin{enumerate}
\item the vertices and edges intersect in disjoint line segments;
\item each such line segment lies on the boundary of precisely one
vertex and precisely one edge;
\item every edge contains exactly two such line segments.
\end{enumerate}
\end{definition}

It is well known and easily seen that ribbon graphs are equivalent to 
cellularly embedded graphs. Intuitively, if $G$ is a cellularly 
embedded graph, a ribbon graph representation results from taking a 
small neighbourhood of $G$.  Neighbourhoods of vertices of $G$ form 
the vertices of the ribbon graph, and neighbourhoods of the edges of 
$G$ form the edges of the ribbon graph. 
On the other hand, if $G$ is a ribbon graph, we simply sew discs into 
each boundary component of the ribbon graph to get a graph cellularly 
embedded in a surface.  Since ribbon 
graphs and cellularly embedded graphs are equivalent we can, and 
will, move freely between them.

Two ribbon graphs are considered to be equal if their corresponding 
embedded graphs are equal. This means that two ribbon graphs are 
equal if there is a  homeomorphisms between their underlying surfaces 
that preserve the vertex-edge structure.  Again, at times we abuse 
notation and identify a ribbon graph with its equivalence class under 
equality.

Just as with graphs, if $G$ is a ribbon graph and $A\subseteq E(G)$, 
then $G-A$ is the ribbon graph obtained from $G$ by deleting all of 
the edges in $A$. Note that $G-A$ is also a ribbon graph and 
therefore describes a cellularly embedded graph. It is this closure 
of the set of ribbon graphs under deletion of edges that makes them 
useful here; note that deleting edges in a cellularly embedded graph 
may result in a non-cellularly embedded graph.
Furthermore,  if $G$ is a ribbon graph and $A\subseteq E(G)$, then  
 $G|_A$ denotes the ribbon subgraph of $G$ induced by $A$, {\em i.e.}, its edge set is $A$ 
and its vertex set consists of all vertices of $G$ which are 
incident to an edge in $A$.  If $G$ is a cellularly embedded graph 
then $G|_A$ is defined to be the cellularly embedded graph 
corresponding to this ribbon graph.

We will need to be able to delete edges from a ribbon graph without 
losing any information about the position of the edge. We will do 
this by recording the position of the edge using labelled arrows.

\begin{definition}
An {\em arrow-marked ribbon graph} consists of a ribbon graph  
equipped with a collection of  coloured arrows, called {\em marking 
arrows}, on the boundaries of its vertices. The marking arrows are 
such that no marking arrow meets an edge of the ribbon graph, and    
there are exactly two marking arrows of each colour.
\end{definition}

\begin{figure}
\begin{center}
\begin{tabular}{ccccc}
\includegraphics[height=15mm]{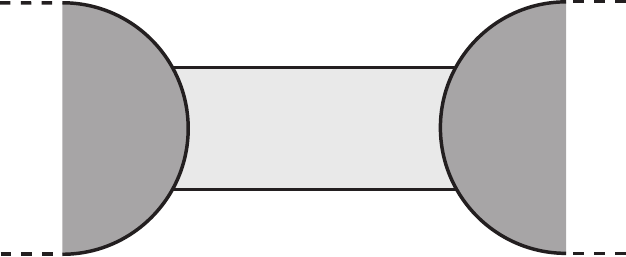}  &
\raisebox{6mm}{\includegraphics[width=11mm]{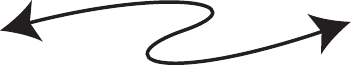}} &
\includegraphics[height=15mm]{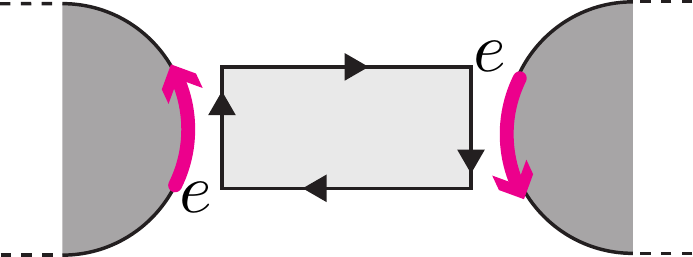}  &
\raisebox{6mm}{\includegraphics[width=11mm]{doublearrow}} &
\includegraphics[height=15mm]{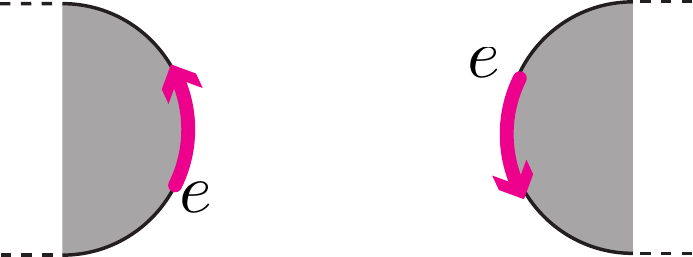} 
\\
$G\,\vec{+}\{e\}$&&&&$G\,\vec{-}\{e\}$
\end{tabular}
\end{center}
\caption{Constructing $G\,\vec{-}\{e\}$ and $G\,\vec{+}\{e\}$. }
\label{arrows}
\end{figure}

Let $G$ be a ribbon graph and $A\subseteq E(G)$. Then we let 
$G\,\vec{-} \, A$ denote the arrow-marked ribbon graph obtained, for 
each edge $e\in A$,  as follows: arbitrarily orient the boundary of 
$e$; place an arrow on each of the two arcs where $e$ meets vertices 
of $G$, such that the directions of these arrows follow the 
orientation of the boundary of $e$; colour the two arrows with $e$; 
and delete the edge $e$. This process is illustrated locally at an 
edge in Figure~\ref{arrows}.

On the other hand, given an arrow-marked ribbon graph $G$ with set of labels 
$A$, we can recover a ribbon graph  $G\vec{+}A$ as follows: for each 
label $e\in A$ take a disc and orient its boundary arbitrarily; add 
this disc to the ribbon graph by choosing two non-intersecting arcs 
on the boundary of the disc and two marking arrows of the same 
colour, and then identifying the arcs with the marking arrows 
according to the orientation of the arrow. The disc that has been 
added forms an edge of a new ribbon graph. Again, this process is 
illustrated in Figure~\ref{arrows}.

\begin{figure}
\centering
\subfigure[A ribbon graph $G$.]{
\includegraphics[scale=.6]{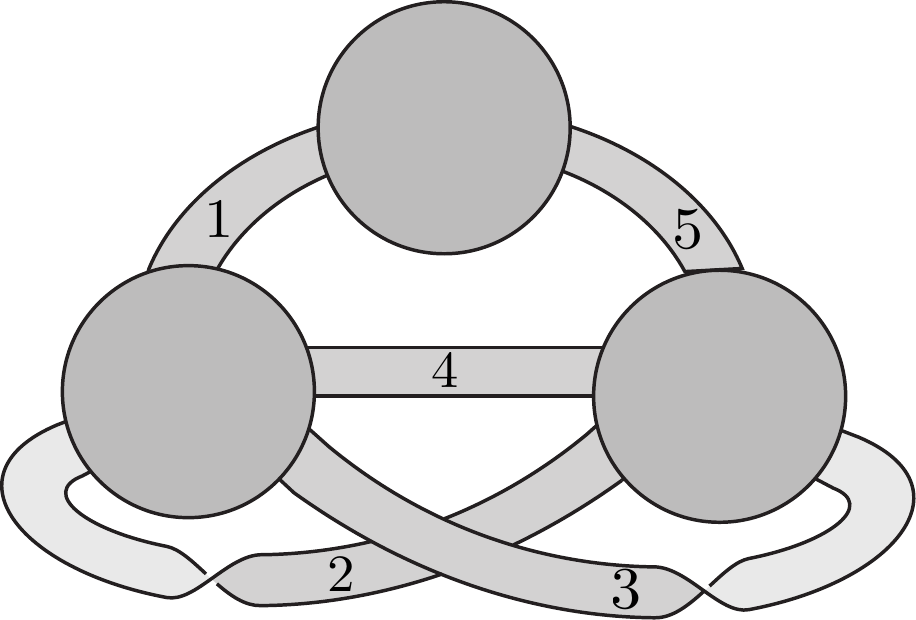}
\label{e.minusb}
}
\hspace{10mm}
\subfigure[$G\,\vec{-}A$ with $A=\{1,2,5\}$  ]{
\includegraphics[scale=.6]{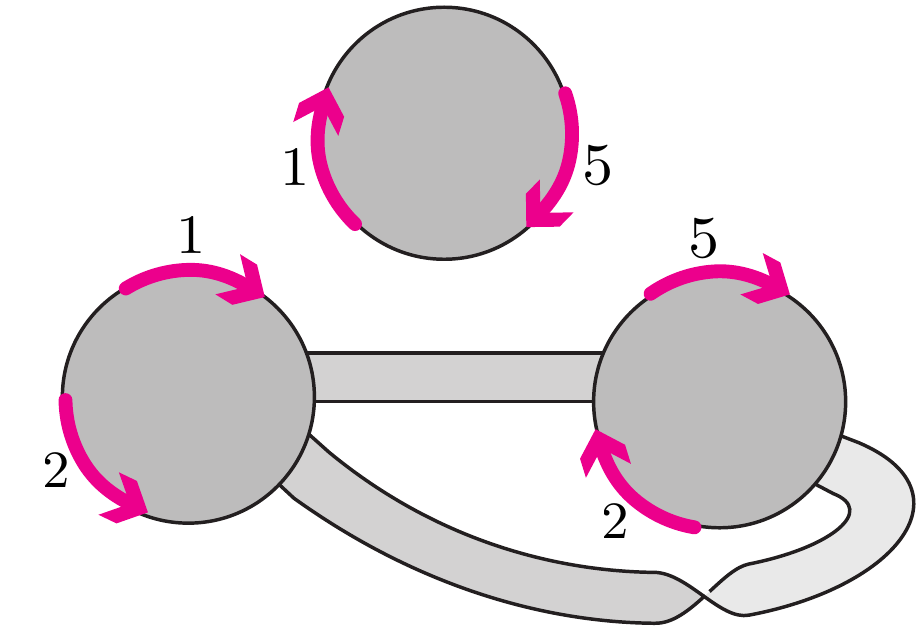}
\label{e.minusc}
}\caption{Two descriptions of the same ribbon graph.}
\label{e.minusa}
\end{figure}

\begin{example}\label{e.minus}
Figure~\ref{e.minusa} shows a ribbon graph $G$ and  its description as the arrow-marked ribbon graph $G\,\vec{-}A$, where $A=\{1,2,5\}$. 
Note that $G$ can be recovered from 
$G\,\vec{-}A$ by taking $A=\{1,2,5\}$ to be the 
set of labels and forming  $(G\,\vec{-}A)\,\vec{+}A$.
\end{example}

From the above we see that every arrow-marked ribbon graph gives rise 
to a ribbon graph. We then say that two arrow-marked ribbon graphs 
are  equal if the ribbon graphs they describe are 
equal.
We will generally abuse notation and regard the set of labels of an 
arrow-marked ribbon graph as a set of edges. This will allow us to 
view $A$ as an edge set in expressions like $G= \left(G\,\vec{-} \,A\right) \vec{+} A$.

\subsection{Geometric duals}

The construction of the {\em geometric dual}, $G^*$, of a cellularly 
embedded graph $G\subset \Sigma$ is well known: $G^*$ is formed  by 
placing one vertex in each face of $G$ and embedding an edge of $G^*$ 
between two vertices whenever the faces of $G$ they lie in are 
adjacent. Observe that $G^*$ has a natural cellular  embedding in 
$\Sigma$, and that there is a natural (cellular) immersion of $G \cup 
G^*$ where each edge of $G$ intersects exactly one edge of $G^*$ at 
exactly one point. We will call this immersion the {\em standard 
immersion} of $G \cup G^*$.

There is a natural bijection between $E(G)$ and $E(G^*)$. We will generally use this bijection to identify the edges 
of $G$ and the edges of $G^*$. However, at times we will be working  
with $G \cup G^*$, so to avoid confusion we will use $e^*$ to denote 
the edge of $G^*$ which corresponds to the edge $e$ of $G$.

Geometric duals have a particularly neat description in the language 
of ribbon graphs. Let $G=(V(G), E(G))$ be a ribbon graph. We can 
regard $G$ as a punctured surface. By filling in the punctures using 
a set of discs denoted $V(G^*)$, we obtain a surface without 
boundary. The {\em geometric dual} of $G$  is the ribbon 
graph $G^* = (V(G^*), E(G))$.

Suppose now that $G$ is an arrow-marked ribbon graph, so that $G$ has 
labelled arrows on its vertices.  Then in the formation of $G^*$ as 
described above, the boundaries of the vertices of $G$ and $G^*$ 
intersect, and  therefore the marking arrows on $G$ induce marking 
arrows on $G^*$. The geometric dual $G^*$ of an arrow-marked ribbon 
graph $G$ is the geometric dual of the underlying ribbon graph 
equipped with the induced marking arrows.

Note that for ribbon graphs geometric duality acts 
disjointly on connected components, so that $(G \sqcup H)^* = G^* 
\sqcup H^*$.

\medskip

We will also need to form geometric duals of non-cellularly embedded 
graphs. Since the properties of duality depend upon whether or not a 
graph is cellularly embedded, we will avoid confusion by denoting the 
dual of a not necessarily cellularly  embedded graph by $G^{\du}$. The embedded 
graph $G^{\du}$ is formed just as the geometric dual of an embedded 
graph is formed but by placing a vertex in each region of $G$, rather 
than each face. That is, if $G\subset \Sigma$ is an embedded graph 
(the embedding may or may not be cellular here), then $G^{\du} 
\subset \Sigma$ is the embedded graph formed by placing   one 
vertex in each region of $G$, and embedding an edge of $G^{\du}$ 
between two vertices whenever the regions of $G$ they lie in are 
adjacent. It is important to note that in general $(G^{\du})^{\du}\neq G$. 
Also, as there are some choices of where to place the edges in its formation, the embedding of $G^{\du}$ is not unique. This fact does not cause any problems here.

\subsection{Partial duals of embedded graphs}\label{ss.pddef}

We can now describe partial duality, which was introduced by 
S.~Chmutov in \cite{Ch1} to unify various results which realize the 
Jones polynomial as a graph polynomial (two of these results were 
first related in \cite{Mo2}). We will use the definition of 
a partial dual  from \cite{Mo4}. 
Chmutov's original (and equivalent) definition of a partial dual can be 
found in \cite{Ch1}.

Let $G$ be a ribbon graph and $A\subseteq E(G)$. Then the partial 
dual $G^A$ of $G$ is formed by: `hiding' the edges that are not in 
$A$ by replacing them with marking arrows using $G\,\vec{-}A^c$; 
forming the geometric dual $(G\,\vec{-}A^c)^*$ (so that the dual is 
only taken with respect to the edges of $G$ that are in $A$); then 
putting back in the edges that are not in $A$, giving 
$(G\,\vec{-}A^c)^*\vec{+} A^c$. The resulting ribbon graph is 
$G^A$. Here and henceforth, $A^c$ denotes the complementary edge set 
$E(G)-A$ of $A$. This process is summarized by the following 
definition.

\begin{definition}
Let $G$ be a ribbon graph and $A\subseteq E(G)$. Then the {\em 
partial dual of $G$ with respect to $A$}, denoted by $G^A$, is given 
by
\[ G^A:= \left(G\,\vec{-} \,A^c\right)^* \vec{+}\, A^c.\]
\end{definition}
 
The partial dual of a cellularly embedded graph is obtained by 
translating into the language of ribbon graphs, forming the partial 
dual, and translating back into the language of cellularly embedded 
graphs.

\begin{figure}
\centering
\subfigure[A ribbon graph $G$.]{
\includegraphics[scale=1.1]{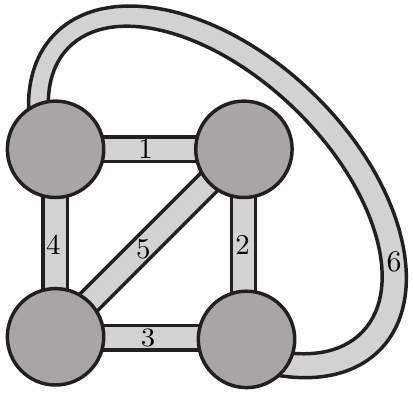}
\label{f.pda}
}
\hspace{10mm}
\subfigure[$G\,\vec{-} \,A^c=$ with $A=\{4,6\}$.  ]{
\includegraphics[scale=1.1]{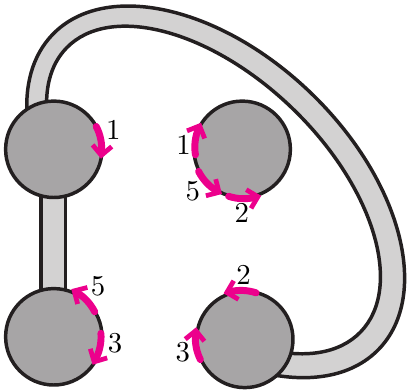}
\label{f.pdb}
}
\hspace{10mm}
\subfigure[$\left(G\,\vec{-} \,A^c\right)^*$.  ]{
\includegraphics[scale=1.1]{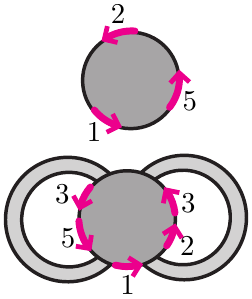}
\label{f.pdc}
}
\subfigure[$\left(G\,\vec{-} \,A^c\right)^* \vec{+}\, A^c=G^A$  ]{
\includegraphics[scale=1.2]{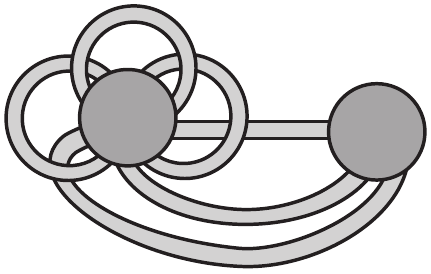}
\label{f.pdd}
}
\caption{Forming a partial dual.}
\label{f.pd}
\end{figure}

\begin{example} Consider the ribbon graph $G$ shown in Figure~\ref{f.pda}. To form the partial dual $G^A$, with $A=\{4,6\}$, first form the arrow-marked ribbon graph $G\,\vec{-} \,A^c$,  as in Figure~\ref{f.pdb}. Then form its geometric dual  $\left(G\,\vec{-} \,A^c\right)^*$,  shown in Figure~\ref{f.pdc}, noting that the labelled arrows on the vertices of $G\,\vec{-} \,A^c$ induce some on $\left(G\,\vec{-} \,A^c\right)^*$. The corresponding ribbon graph $\left(G\,\vec{-} \,A^c\right)^* \vec{+}\, A^c$ is the partial dual $G^A$ and is shown in Figure~\ref{f.pdd}.
\end{example}

Further examples of partial duals can be found in \cite{Ch1, Mo3,Mo4, 
Mo5}, and the references therein.

We will need the following basic properties of partial duality later. These properties are due to Chmutov and can be found in~\cite{Ch1}.
\begin{proposition}\label{p.pd2}
Let $G$ be a  ribbon graph and $A, B\subseteq E(G)$.  Then 
\begin{enumerate}
\item $G^{\emptyset}=G$;
\item  $G^{E(G)}=G^*$, where $G^*$ is the geometric dual of $G$;
\item $(G^A)^B=G^{A\Delta B}$, where $A\Delta B := (A\cup 
B)\backslash (A\cap B)$ is the symmetric difference of $A$ and $B$.
\end{enumerate}
\end{proposition}

\section{Eulerian and Bipartite partial duals}\label{s.bepd}

This section gives our first main result, which appears as Theorem~\ref{t.pde} (the second main result being Theorem~\ref{t.bipmed}). This is the extension to partial duality of the classical result stated in Theorem~\ref{t.intro1}. The relationship in Theorem~\ref{t.intro1} between bipartite and Eulerian graphs is usually given for connected plane graphs, but here we need  the following slightly more general 
form.
\begin{theorem}\label{t.euler}
Let $G$ be a graph embedded in the plane. Then 
\begin{enumerate}
\item the components of $G$ are all Eulerian if 
and only if $G^{\du}$ is bipartite; and
\item  $G$ is bipartite if 
and only if $G^{\du}$ is Eulerian.
\end{enumerate}
\end{theorem}
\begin{proof}
Let $G_1, \cdots , G_k$ denote the components of the embedded graph $G$. Then the plane
embedding of $G$ can be obtained by forming the connected sum of cellular plane embeddings of $G_1, \cdots , G_k$. (The duals $G_i^*$ below are formed with respect the these embeddings.)  In terms of duals, this means that $G^{\du}$ can be obtained by amalgamating $G_1^*, \cdots , G_k^*$ at vertices (two vertices of the duals $G_i^*$ and $G_j^*$ being amalgamated if, in the construction of $G$, a connected sum involves the corresponding faces of the plane graphs $G_i$ and $G_j$). If the components of $G$ are Eulerian, then so are $G_1, \cdots , G_k$. By Theorem~\ref{t.intro1}, it follows that $G_1^*, \cdots , G_k^*$ are bipartite, and since amalgamating bipartite graphs at a vertex results in a bipartite graph, $G^{\du}$ is bipartite. Conversely, if $G^{\du}$ is bipartite, then so are $G_1^*, \cdots , G_k^*$. By Theorem~\ref{t.intro1} it follows that $G_1, \cdots , G_k$ are Eulerian, and therefore so are the components of $G$.

The second item in the theorem follows by interchanging the words bipartite and Eulerian in the above argument.
\end{proof}

We begin with the observation that biparticity is a property of 
abstract graphs rather than embedded graphs. Accordingly, in order to 
study the $2$-colourability of partially dual embedded graphs, it 
suffices to study their underlying abstract graphs. This allows us to 
use tools developed in \cite{Mo4} for partial duals of abstract 
graphs to prove Theorem~\ref{t.pde}. Also, to prove the theorem, we introduce a new way of 
constructing abstract graphs that are partial duals.

Recall that an embedded graph $G$ consists of an embedding of a graph 
$\hat{G}$ into a surface. We call the graph $\hat{G}$ the {\em 
underlying abstract graph} of $G$. If $G$ and $H$ are two embedded 
graphs with the same underlying abstract graphs we will say that $G$ 
and $H$ are {\em equivalent as abstract graphs} and write $G\cong H$. 
The notion of partially dual abstract graphs was introduced in 
\cite{Mo4}.

\begin{definition}
Two abstract graphs are said to be {\em partial duals} if they are 
the underlying abstract graphs of two partially dual embedded 
graphs.
\end{definition}

\begin{remark}
It is important to observe that although partial duality is a 
transitive relation for embedded graphs, it is not a transitive 
relation for abstract graphs. For example, the two abstract graphs 
\raisebox{-3mm}{\includegraphics[height=8mm]{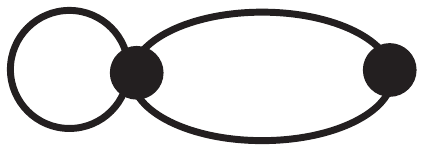}}
and
\raisebox{-6mm}{\includegraphics[height=15mm]{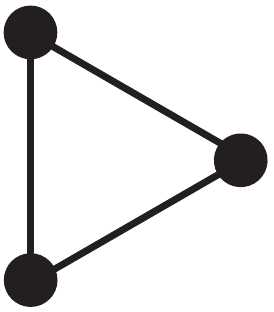}}
are partial duals because they are the underlying abstract graphs of 
the partial duals
\raisebox{-5mm}{\includegraphics[height=15mm]{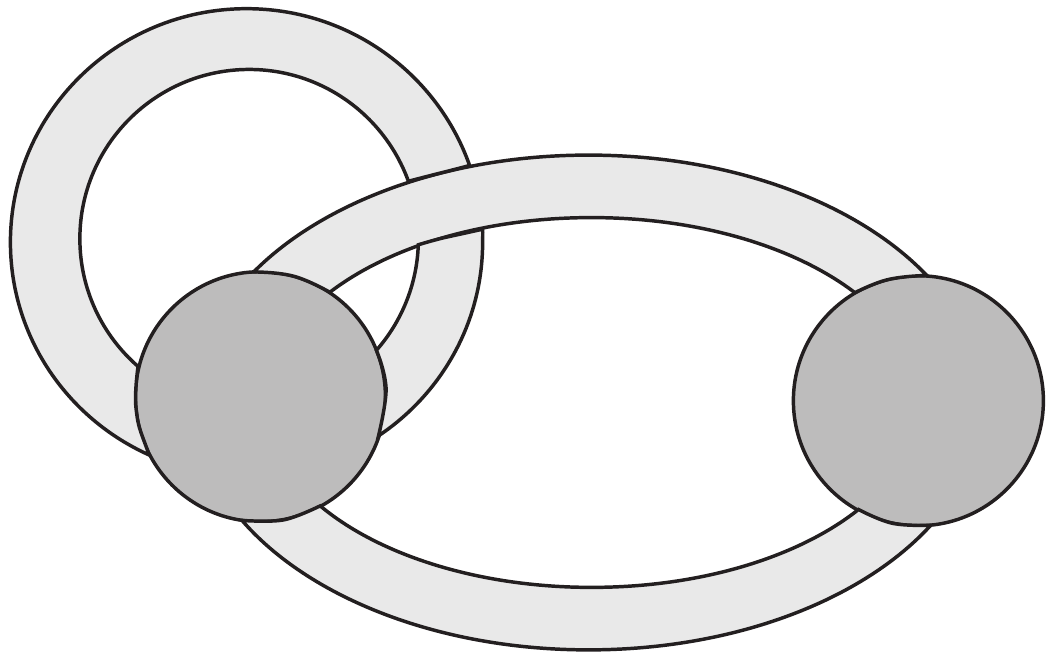}}
and
\raisebox{-6mm}{\includegraphics[height=15mm]{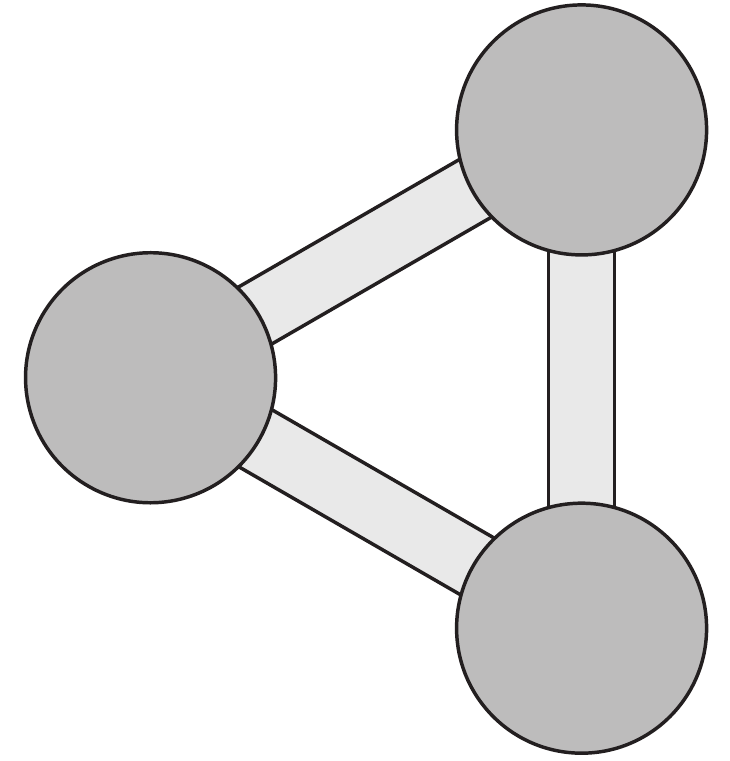}} respectively.
Also, the abstract graphs 
\raisebox{-3mm}{\includegraphics[height=8mm]{pdem2}}
and
\raisebox{-3mm}{\includegraphics[height=8mm]{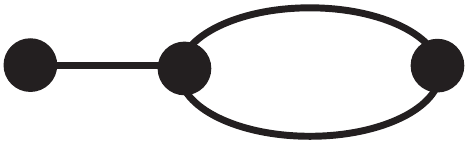}}
are partial duals, since they are the underlying abstract graphs of
\raisebox{-5mm}{\includegraphics[height=12mm]{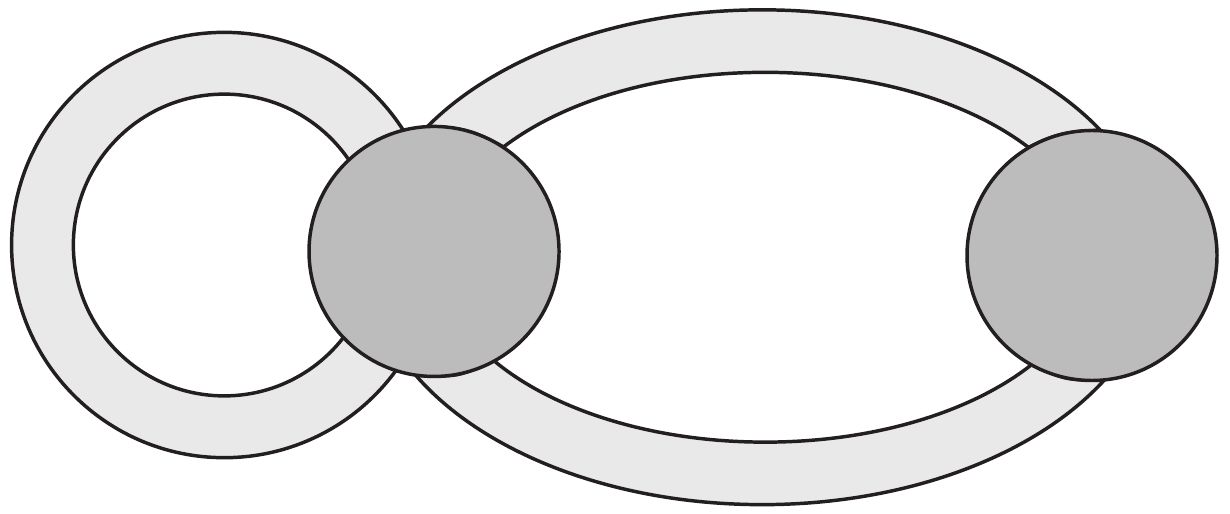}}
and
\raisebox{-5mm}{\includegraphics[height=12mm]{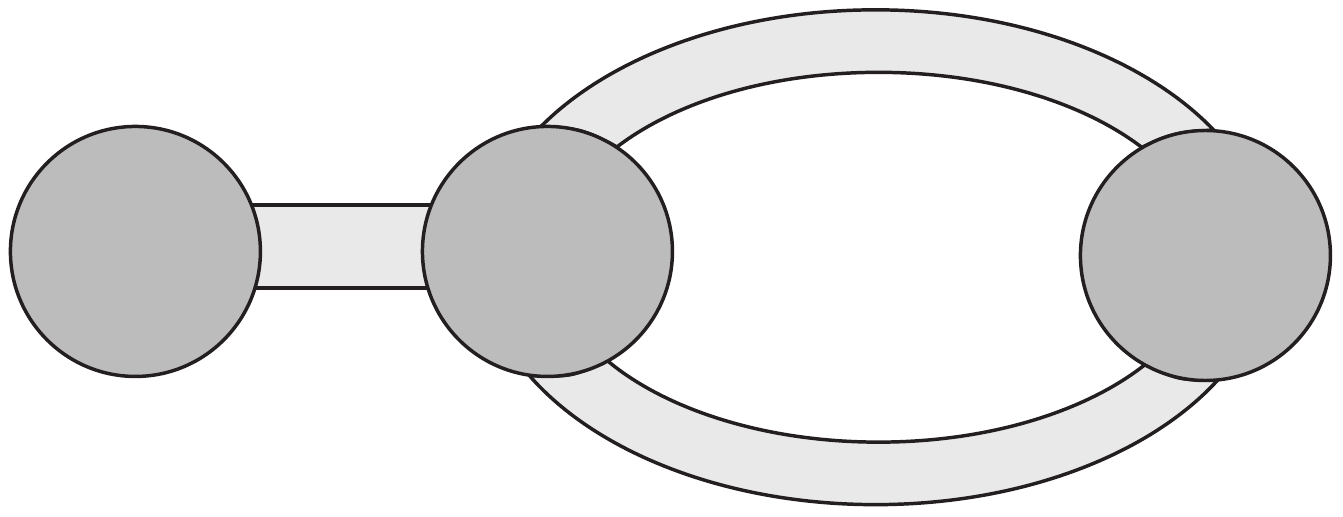}} respectively.
However, 
\raisebox{-6mm}{\includegraphics[height=15mm]{pdem3}}
and 
\raisebox{-3mm}{\includegraphics[height=8mm]{pdg1}}
are not partially dual abstract graphs. This observation has implications for the results presented here.
\end{remark}

We now give a new way of constructing partially dual abstract graphs. 
Theorem~\ref{t.pde} will follow easily from this construction. Given 
an embedded graph $G$, Theorem~\ref{t.newpds} provides a  way to 
obtain an embedded graph that is equivalent to $G^A$ as an abstract 
graph but not necessarily as an embedded graph. This result is 
especially useful here since $G^A$ is bipartite if and only if any 
embedded graph $H$ with $H\cong G^A$ is bipartite, so we need not worry about the  embedding of $G^A$.

Since we will be working simultaneously with an embedded graph $G$ 
and its dual $G^*$, we will use  a superscript `$\ast$' to denote 
corresponding edges and edge sets in $G^*$. For example, if $e$ is 
and edge in $G$ then $e^*$ denotes the edge in $G^*$ under the 
natural identification of $E(G)$ with $E(G^*)$.

\begin{theorem}\label{t.newpds}
Let $G$ be a connected, cellularly embedded graph, and $A\subseteq 
E(G)$. Then
\[ [ (G\cup G^{*}) - (A^c  \cup  A^{*} )  ]^{\du}  \cong G^A.  \]
Here,  $G\cup G^{*}$ has the standard immersion. 
\end{theorem}

In order not to disrupt our narrative on bipartite partial duals, we 
defer the somewhat technical proof of Theorem~\ref{t.newpds} until 
Section~\ref{s.proof}.

We emphasize the fact that in general $[(G\cup G^{*}) - (A^c \cup 
A^{*})]^{\du}$ and $G^A$ are not equal as embedded graphs. This is 
what makes the above theorem significant: we have found a way of 
constructing partial duals of abstract graphs that does not require 
us to pass through partially dual ribbon graphs. It is perhaps 
prudent to highlight a second point, that in general $([(G\cup G^{*}) 
- (A^c \cup A^{*})]^{\du})^*$ and $(G^A)^*$ are not isomorphic as 
abstract graphs. This means that, for a plane graph $G$,  if $[(G\cup G^{*}) - (A^c \cup 
A^{*})]^{\du} \cong G^A$ is bipartite (respectively Eulerian), then 
although the dual $([(G\cup G^{*}) - (A^c \cup A^{*})]^{\du})^*$ is 
Eulerian (respectively bipartite), we do not know whether $(G^A)^*$ 
is bipartite or Eulerian.

\begin{figure}
\centering
\subfigure[A plane graph $G$.]{
\includegraphics[scale=0.6]{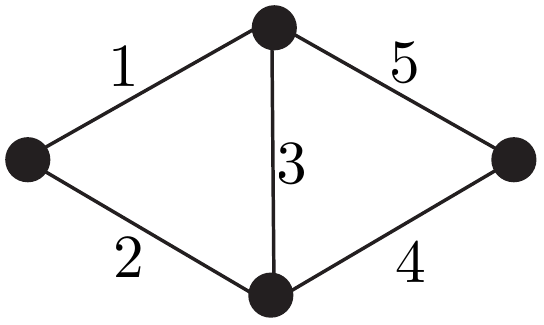}
\label{f.nnpda}
}
\hspace{10mm}
\subfigure[$G\cup G^{*}$.  ]{
\includegraphics[scale=0.6]{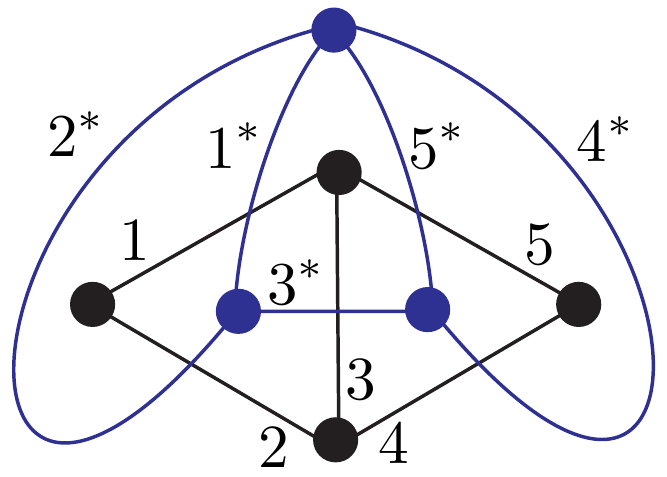}
\label{f.nnpdb}
}
\hspace{10mm}
\subfigure[$(G\cup G^{*}) - (A^c \cup A^{*})$.  ]{
\includegraphics[scale=.6]{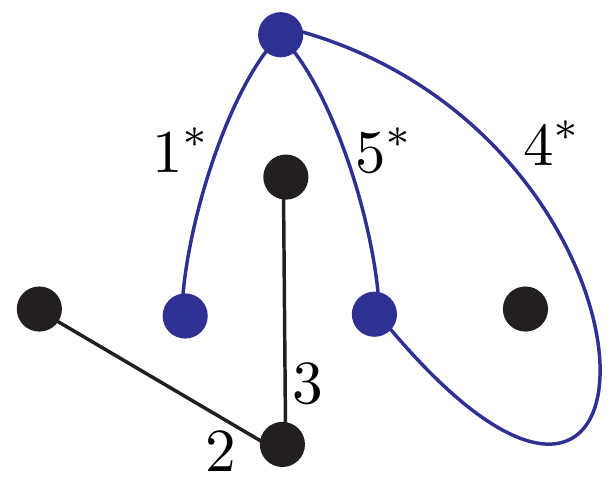}
\label{f.nnpdc}
}
\subfigure[Forming the geometric dual.]{
\includegraphics[scale=.6]{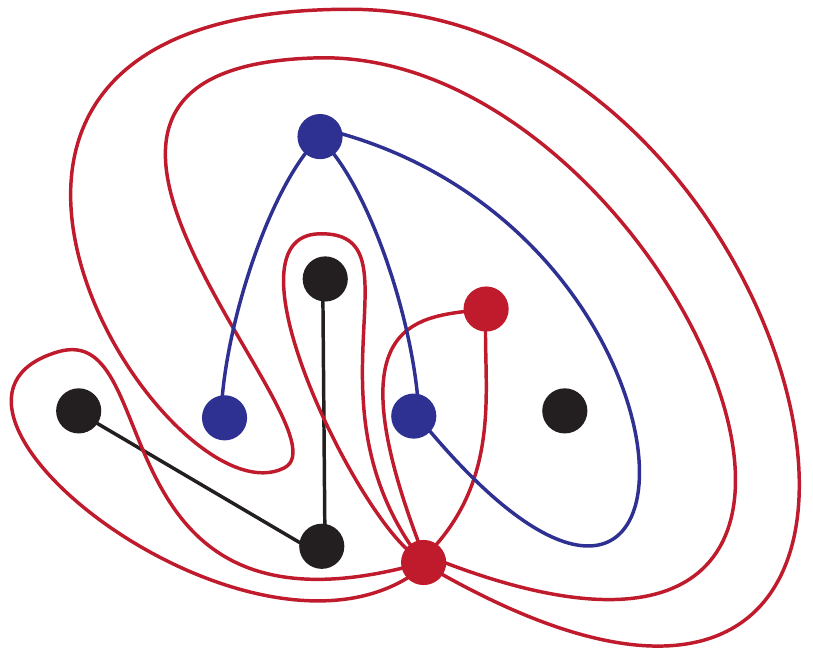}
\label{f.nnpdd}
}
\hspace{10mm}
\subfigure[$\left((G\cup G^*) - (A^c \cup A^*)\right)^{\du}$.]{
\includegraphics[scale=0.6]{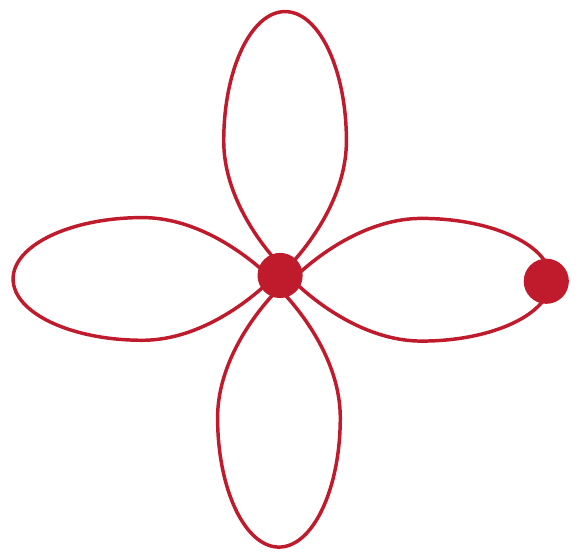}
\label{f.nnpde}
}
\hspace{10mm}
\subfigure[The partial dual $G^A$ of $G$. ]{
\includegraphics[scale=1.1]{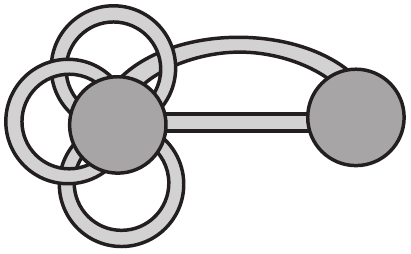}
\label{f.nnpdf}
}
\caption{Forming partially dual abstract graphs using Theorem~\ref{t.newpds}.}
\label{f.nnpd}
\end{figure}

\begin{example}
Theorem~\ref{t.newpds} is illustrated in Figure~\ref{f.nnpd}. For the  cellularly embedded graph $G$  shown  in Figure~\ref{f.nnpda}, $G\cup G^{*}$ is shown in Figure~\ref{f.nnpdb}. 
 Taking  $A=\{2,3\}$, we have $A^c=\{1,4,5\}$ and $A^{*}=\{2^{*},3^{*}\}$. With these sets, $(G\cup G^{*}) - (A^c \cup A^{*})$ is shown in Figure~\ref{f.nnpdc}. Figure~\ref{f.nnpde} illustrates the formation of the geometric dual  $[(G\cup G^{*}) - (A^c \cup A^{*})]^{\du}$, which is given in Figure~\ref{f.nnpde}. 
 
On the other hand, by regarding $G$ as a genus zero ribbon graph, we 
can form the partial dual $G^A$  to obtain the ribbon graph in Figure~\ref{f.nnpdf}, which is equivalent as an abstract graph, but not as an embedded graph, to the graph $[(G\cup G^{*}) - (A^c \cup A^{*})]^{\du}$ in Figure~\ref{f.nnpde}. 
\end{example}

Our generalization of Theorem~\ref{t.euler} is the following theorem, 
in which we identify the edges of $G$ and $G^A$ in the standard way. 
\begin{theorem}\label{t.pde}
Let $G$ be a plane graph and $A\subseteq E(G)$. Then
\begin{enumerate}
\item \label{te1} $G^A$ is bipartite if and only if the components 
of $G|_{A}$  and $G^*|_{A^c}$ are Eulerian;
\item \label{te2} $G^A$ is Eulerian if and only if $G|_{A}$ and 
$G^*|_{A^c}$ are bipartite.
\end{enumerate}
\end{theorem}

Note that Theorem~\ref{t.euler} is obtained from Theorem~\ref{t.pde} 
by setting $A =\emptyset$. Also note that, as one would expect, 
Theorem~\ref{t.pde} does not hold for graphs embedded in higher genus 
surfaces.

\begin{proof}[Proof of Theorem~\ref{t.pde}]
Let $G$ be a plane graph, and set $\Phi:= (G\cup G^{*}) - (A^c \cup 
A^{*})$. By Theorem~\ref{t.newpds}, $G^A$ is bipartite if and only if 
 $\Phi^{\du}$ is bipartite. But, as $\Phi$ is embedded 
in the plane, Theorem~\ref{t.euler} implies that  
$\Phi^{\du}$ is bipartite if and only if each component of  $\Phi$ is 
Eulerian. This happens if and only if each component of $G-A^c$ and 
of $G^*-A^*$ is Eulerian, which happens if and only if each component 
of $G|_A$ and of $G^*|_{A^c}$ is Eulerian. This argument can be 
restated with the words bipartite and Eulerian exchanged, thus 
completing the proof.
\end{proof}

\begin{remark}
The obvious extension to Theorem~\ref{t.euler} is that, for a plane 
graph $G$, $G^A$ is bipartite if and only if $(G^A)^*$ is Eulerian. 
However, this statement is not true in one direction. For example, if 
$G^A$ is the cellularly embedded graph on a torus that consists of a 
loop added to a $2$-cycle, then $G^A$ is Eulerian, but $(G^A)^*$ is 
not bipartite. The statement is true in the other direction since the 
dual of any bipartite embedded graph is Eulerian (see 
Remark~\ref{r.btoe}).
\end{remark}

\section{Bipartite graphs and circuits in medial graphs}\label{s.med}

In this section we apply Theorem~\ref{t.pde} to obtain our second 
main result, which is a characterization of those edge sets of plane 
graphs which give rise to bipartite partial duals. This 
characterization will be in terms of circuits in the medial graph. 

\subsection{Medial graphs}

The {\em medial graph} $G_m$ of a plane graph $G$ is the $4$-regular 
plane graph obtained from $G$ by placing a vertex on each edge of 
$G$, and joining two such vertices by an edge embedded in a face 
whenever the two edges on which they lie are on adjacent edges of the 
face.

\begin{figure}
    \subfigure[A plane graph $G$.] 
    {
    \label{f.m1}\raisebox{5mm}{\includegraphics[height=3cm]{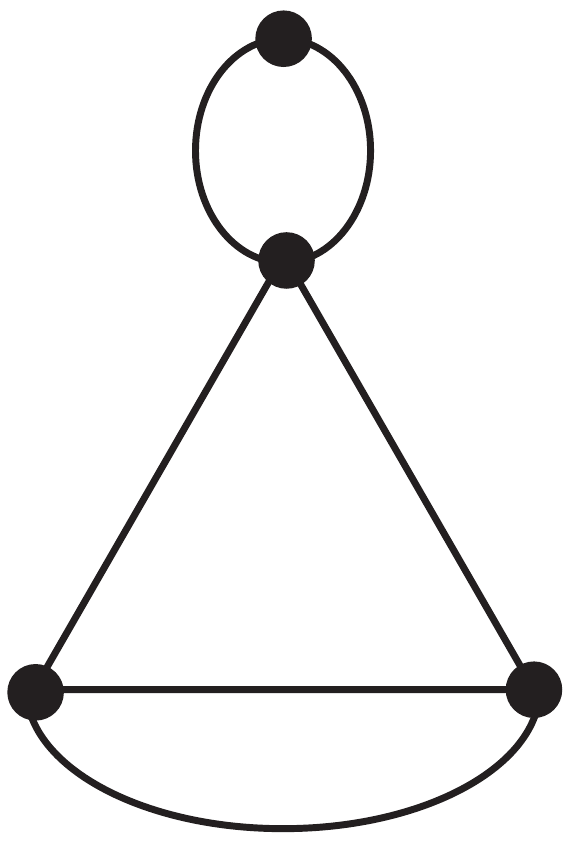}}\hspace{5mm}
    }
    \hspace{5mm}
    \subfigure[Its canonically face $2$-coloured medial graph $G_m$.]
    {
    \label{f.m2}\includegraphics[height=4cm]{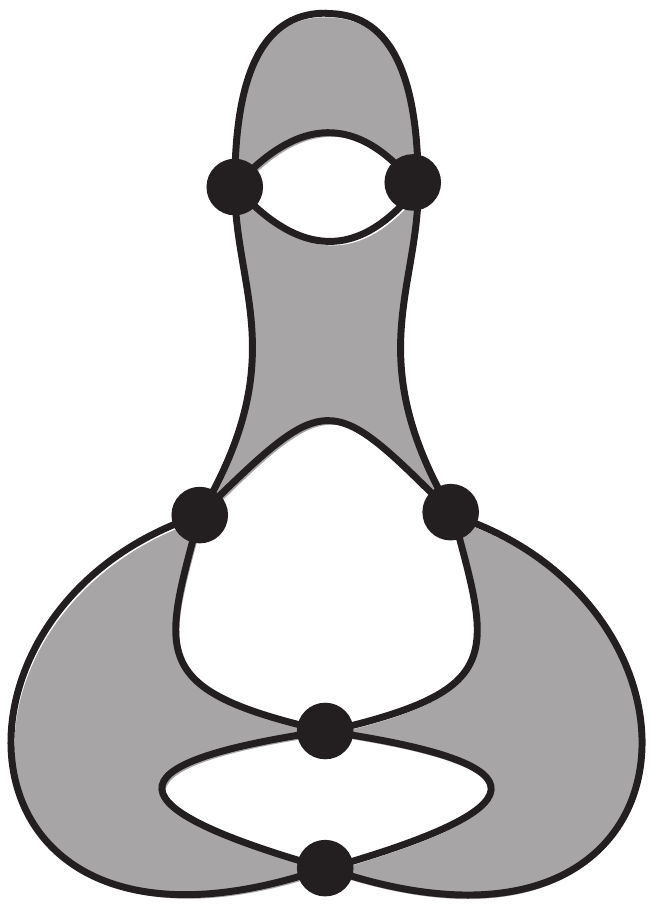}\hspace{5mm}
    }
    \hspace{5mm}
    \subfigure[An all-crossing direction of $G_m$. Here $c(G_m)=1$.]
    {
    \label{f.m3}\includegraphics[height=4cm]{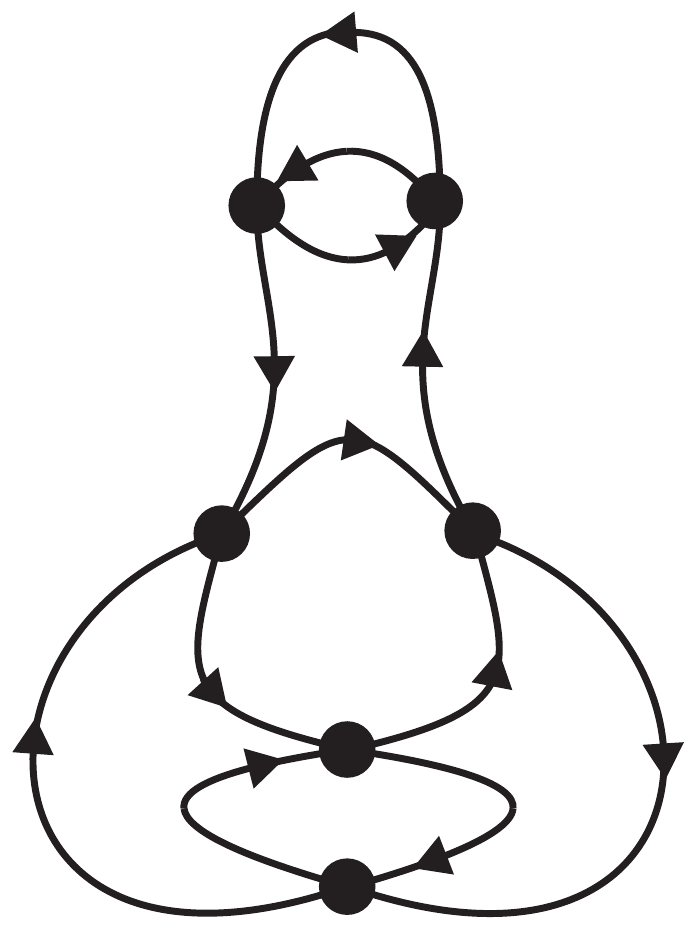}\hspace{5mm}
    }
    \hspace{5mm}
    \subfigure[The corresponding $\{c,d\}$-labelling on $G$.]
    {
    \label{f.m4}\raisebox{5mm}{\includegraphics[height=3cm]{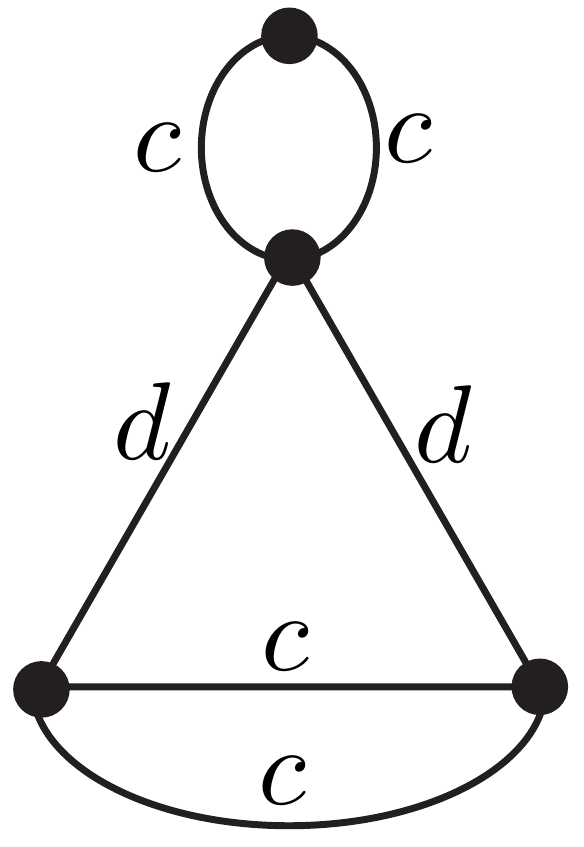}}\hspace{10mm}
    }
\caption{Examples of constructions associated with medial graphs.}
\label{f.m}
\end{figure}

Each vertex of $G$ corresponds to a face of $G_m$. If we colour all 
such faces of $G_m$ black and the remaining faces white we obtain a 
face $2$-colouring of $G_m$ which we will call the {\em canonical 
face $2$-colouring}. (See Figures \ref{f.m1} and \ref{f.m2}.)

We are interested in particular directed graphs which arise by 
directing the edges of a medial graph. An {\em all-crossing 
direction} of $G_m$ is an assignment of a direction to each edge of 
$G_m$ in such a way that at each vertex $v$ of $G_m$, when we follow 
the cyclic order of the directed edges incident to $v$, we find 
(head, head, tail, tail). (See Figure~\ref{f.m3}.)

We will let $c(G_m)$ denote the number of circuits in any 
all-crossing direction of $G_m$ which are obtained by following the 
directed edges in such a way that at each vertex, we enter and exit 
at a head and tail which are not adjacent in the cyclic order of the 
incident edges at that vertex. (See Figure~\ref{f.m3}.) We observe that $c(G_m)$ is 
independent of the choice of all-crossing direction of $G_m$, and 
that $G_m$ admits $2^{c(G_m)}$ all-crossing directions.

\begin{figure}
\begin{center}
\begin{tabular}{ccc}
\includegraphics[height=20mm]{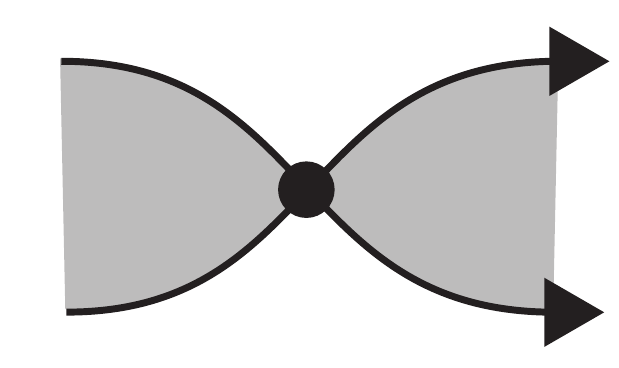} & \hspace{2cm} & 
\includegraphics[height=20mm]{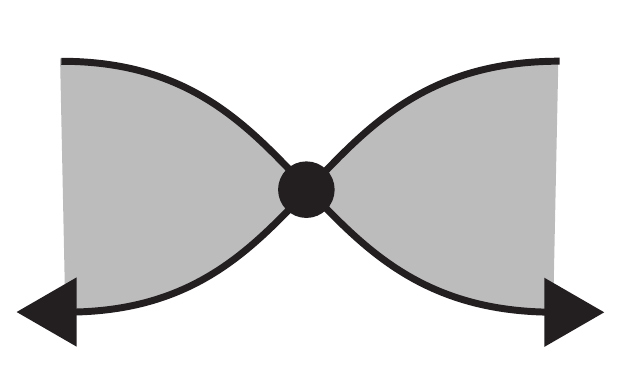} \\ $c$-vertex & & $d$-vertex
\end{tabular}
\end{center}
\caption{The definition of $c$-vertices and $d$-vertices.}
\label{c d vertices}
\end{figure}

If $G_m$ is equipped with the canonical face $2$-colouring then we 
can partition the vertices  of $G_m$ by calling each vertex a {\em 
$c$-vertex} or a {\em $d$-vertex} according to the scheme shown in 
Figure~\ref{c d vertices}.
Furthermore, since the vertices of $G_m$ correspond to edges of $G$, 
each all-crossing direction of $G_m$ gives rise to a 
$\{c,d\}$-labelling of the edges of $G$. We call the edges of $G$ 
which correspond to $c$-vertices of $G_m$ {\em $c$-edges}, and we 
call the edges of $G$ which correspond to $d$-vertices of $G_m$ {\em 
$d$-edges}. (See Figure~\ref{f.m4}.)

We will need the following observation. 
\begin{proposition}\label{p.cd}
Let $G$ be a plane graph, then $e$ is a $c$-edge in $G$ if and only 
if $e^*$ is a $d$-edge in $G^*$.
\end{proposition} 
\begin{proof}
The result follows by observing that the medial graphs $G_m$ and 
$(G^*)_m$ are equal and that the canonical face $2$-colouring of 
$(G^*)_m$ is obtained from that of $G_m$ by switching the colour of 
each face.
\end{proof}

\subsection{Bipartite partial duals and medial graphs}

We now state and prove our second main result, which is a characterization of bipartite partial 
duals in terms of medial graphs.
\begin{theorem}\label{t.bipmed}
Let $G$ be a plane graph. Then the partial dual $G^A$ is bipartite if 
and only if $A$ is the set of $c$-edges arising from an all-crossing 
direction of $G_m$. 
\end{theorem}
We will deduce the theorem from the following lemmas.

\begin{figure}
\includegraphics[height=4cm]{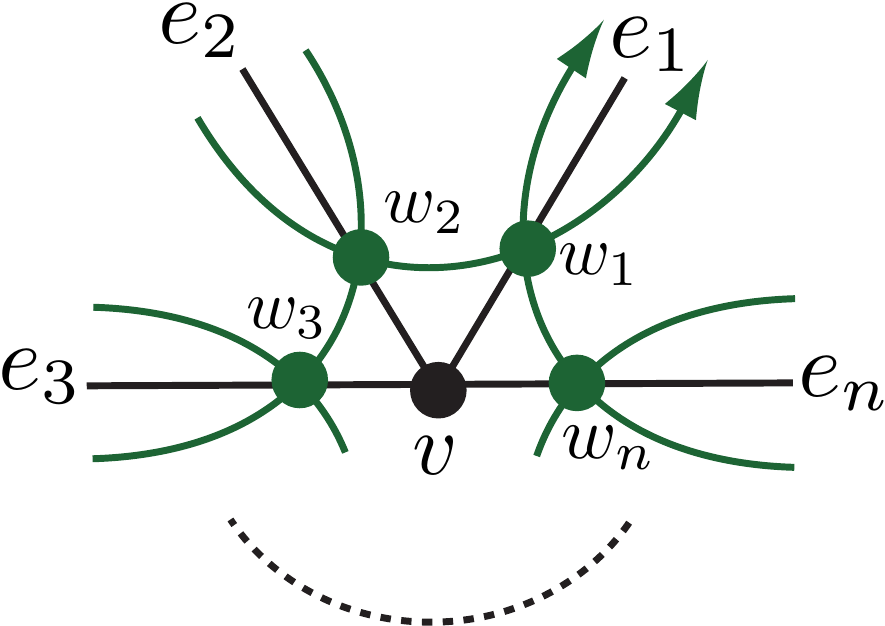}

\caption{Walking around vertex $v$ of graph $G$.}
\label{walk}
\end{figure}

\begin{lemma}\label{l.bm1}
Let $G$ be a plane graph and let $C$ be a set of $c$-edges of $G$ 
arising from an all-crossing direction of $G_m$. Then each component 
of $G|_C$ and of $G^*|_{C^c}$ is Eulerian.
\end{lemma}
\begin{proof}
Consider a vertex $v$ of $G$, with half-edges $e_{1},e_{2},\ldots,e_{n}$ 
adjacent to it. This vertex $v$ corresponds to a face $f$ of $G_{m}$, 
bounded by edges of $G_{m}$ joined by (not necessarily distinct) vertices 
$w_{1},w_{2},\ldots,w_{n}$ corresponding to the half-edges 
$e_{1},e_{2},\ldots,e_{n}$ of $G$. Let us take a position on the 
edge $w_{n}w_{1}$ of $G_{m}$, and then walk around the face $f$ 
until we return to our starting point. We note that because $G_{m}$ 
has an all-crossing direction, each edge will be directed, and 
we start our walk in this given direction.

We walk until we meet our first vertex, which without loss of 
generality we may take to be $w_{1}$. This vertex has a label, either 
$c$ or $d$. Let us assume that it is a $c$-vertex; see Figure 
\ref{walk}. When we cross $w_{1}$ and continue walking on to the next 
half-edge $w_{1}w_{2}$ we will walk against the given direction of 
$w_{1}w_{2}$. In order for us to get back to where we started, the 
edges along which we walk will have to change direction in total 
an even number of times. Therefore there must be an even number of 
$c$-vertices around the face $f$, and so $v$ is adjacent to an even 
number of $c$-half-edges.

Now let us assume that the first vertex we meet $w_{1}$ is a 
$d$-vertex. When we cross this vertex the next edge of $G_{m}$ will 
be directed compatibly with the direction in which we are walking. So 
we keep walking until we meet a $c$-vertex, and then the previous 
argument applies. If we do not encounter any $c$-vertices then all 
the half-edges adjacent to the vertex $v$ are labelled $d$.

Since a graph is Eulerian if and only if each of its vertices is of 
even degree, it follows that each component of $G|_C$ is Eulerian.
Also, it follows from Proposition \ref{p.cd} that each component of 
$G^*|_{C^c}$ is Eulerian.
\end{proof}

(In fact it is not difficult to show the equivalent property that any 
cycle of $G$ contains an even number of $d$-edges.)

\begin{lemma}\label{l.bm2}
Given a $\{c,d\}$-colouring of the edges of a  plane graph $G$, with 
the set of $c$-edges denoted $A$, if each component of $G|_{A}$ is 
Eulerian and each component of $G^{*}|_{A^{c}}$ is Eulerian, then the 
$\{c,d\}$-colouring arises from an all-crossing direction of $G_m$.
\end{lemma}

\begin{figure}

 \subfigure[Extending an orientation from $e^{\du}$.]
    {
      \label{f.w4}\raisebox{2mm}{\includegraphics[height=1.5cm]{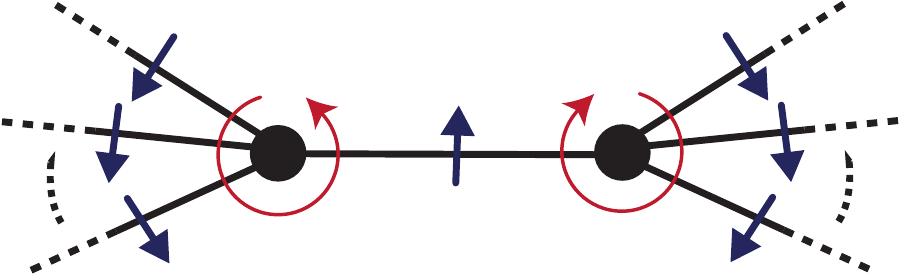} \hspace{5mm}}
    }
    \subfigure[Orientation of neighbouring blocks.]
    {
      \label{f.w2}\includegraphics[height=2cm]{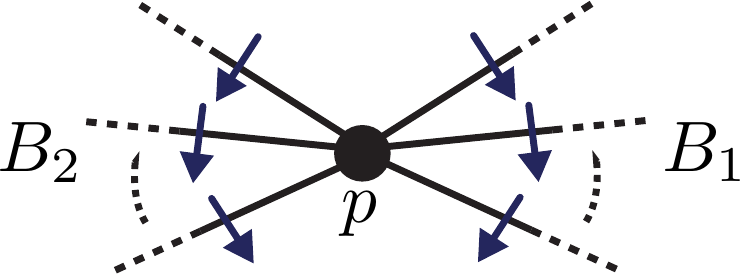}\hspace{5mm}
    } 
    
\label{blocks}
\caption{Extending an orientation of $\Phi^{\du}$.}
\end{figure}

\begin{figure}

\includegraphics[height=1.5cm]{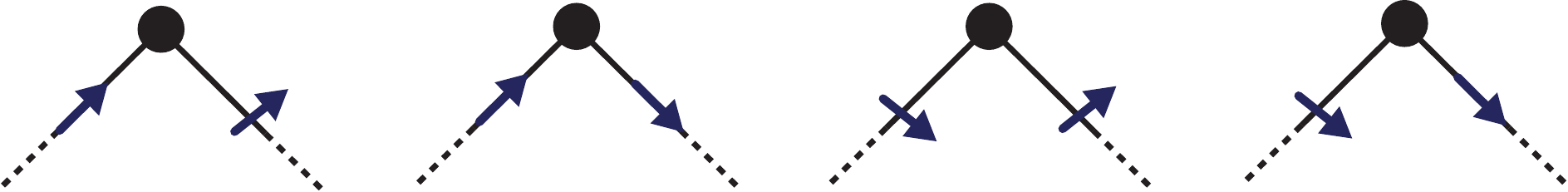}

\caption{Four possible situations for two edges which share both a 
vertex and a region.}
\label{4sits}
\end{figure}

\begin{proof}
The graph $G$ is plane, so it is connected and equipped with an embedding 
$i:G\rightarrow S^{2}$. This induces an embedding
$G^{*}\rightarrow S^{2}$, which we also denote $i$ (so that $i(G)\cup 
i(G^*)$ is the standard immersion). Now take $i(G|_{A})\cup 
i(G^{*}|_{A^{c}})$, and denote the resulting graph by $\Phi$.

We first show that the regions of $\Phi$ must be either discs or 
annuli, by arguing that no two components of $G^{*}|_{A^{c}}$ can lie 
in the same region of $i(G|_{A})$.

Consider a region $f$ of $i(G|_{A})$, drawn in $G$. If it is a region 
of $G$, then it is a disc. If it is not a region of $G$, then it 
contains other edges and vertices of $G$, all the edges being marked 
$d$. These edges and vertices divide $f$ into regions of $G$. But $G$ 
is connected, so we can walk between any two of these regions via 
other such regions. Therefore, the part of $G^{*}|_{A^{c}}$ lying in 
$f$ must be connected, as required.

In $\Phi^{\du}$, this means that ay any separating vertex of $\Phi^{\du}$  exactly two blocks meet.

We can also see that $\Phi^{\du}$ is bipartite, as follows. Each 
component of the (disjoint) graphs $G|_{A}$ and $G^{*}|_{A^{c}}$ is 
Eulerian, and so each vertex of $\Phi$ has even degree. Hence each 
region of $\Phi^{\du}$ has an even number of edges. Any cycle in 
$\Phi^{\du}$ can be formed by adding boundaries of regions mod $2$. 
Hence the result.

Note next that because we will be switching between directed graphs 
and their duals it is natural to consider not only the usual 
``longitudinal'' direction along an edge, but also a ``transverse'' 
direction.

Now choose a block $B_{1}$ of $\Phi^{\du}$. It is bipartite, and has no 
cut vertices. Choose an arbitrary transverse direction on an edge 
$e^{\du}$ in $B_{1}$, and use this to determine a clockwise and 
anti-clockwise orientation of the vertex at each of the ends of $e^{\du}$. (See Figure~\ref{f.w4}.)  Any other 
vertex in $B_{1}$ will be clockwise if it is an even distance from a 
clockwise vertex, and anti-clockwise otherwise. This is consistent 
because $B_{1}$ is bipartite, and so it has no odd cycles.

Next, back in $\Phi^{\du}$, let $p$ be a separating vertex of $B_{1}$. (If 
$B_{1}$ has no separating vertex, then $B_{1}=\Phi^{\du}$ and we skip this 
step.) By our observation above, there is a next block: call it 
$B_{2}$ and orient it as in Figure~\ref{f.w2}. Proceed in this way 
until all the edges of $\Phi^{\du}$ have been given a direction.

Finally, transfer this direction back to the edges of $\Phi$. This 
in turn induces a direction for the edges of $G$, mixed, in the 
sense that the $d$ edges are directed transversely.

The vertices of the medial graph $G_{m}$ inherit their $c$ or $d$ 
status from the edges of $G$, and around each vertex the edges have 
now been directed as in Figure~\ref{c d vertices}. We have to check 
that as we move from one vertex of $G_{m}$ to another along an edge, 
the directions are consistent. To see this, note that for any two 
edges in $G$ which share both a vertex and a region, our direction 
must give one of the four situations in Figure \ref{4sits}, or their 
opposite orientations. So the local directions around each vertex 
of $G_{m}$ do arise from a global direction of the edges of 
$G_{m}$, and we thus have an all-crossing direction, as required.
\end{proof}

\begin{proof}[Proof of Theorem~\ref{t.bipmed}.]
By combining Lemmas~\ref{l.bm1} and \ref{l.bm2}, we have that each 
component of $G|_{A}$ and of $G^{*}|_{A^{c}}$ is Eulerian if and only 
if $A$ is the set $c$-edges of $G$ arising from an all-crossing 
direction of $G_m$. The result then follows by Theorem~\ref{t.pde}.
\end{proof}

\begin{corollary}
Let $G$ be a plane graph. Then $G$ has at most $2^{c(G_m)-1}$ 
bipartite partial duals.  
\end{corollary}
\begin{proof}
$G_m$ admits $2^{c(G_m)}$ all-crossing directions, each direction 
giving rise to a bipartite partial dual by Theorem~\ref{t.bipmed}. 
However, reversing the direction of each edge in an all-crossing 
direction of $G_m$ does not change the $\{c,d\}$-colouring of $G$, 
accounting for the `$-1$' in the exponent. There may be a further 
reduction in the number of bipartite partial duals, as each 
$\{c,d\}$-colouring of $G$ need not result in a distinct partial dual 
of $G$.
\end{proof}

The following corollary of Theorem~\ref{t.bipmed} provides a way of 
constructing some, but not all, of the Eulerian partial duals of a 
plane graph. 
\begin{corollary}\label{c.deu}
Let $G$ be a plane graph. If $A$ is the set of $d$-edges arising from 
an all-crossing direction of $G_m$, then $G^A$ is Eulerian.
\end{corollary}
\begin{proof}
Since $A$ is the set of $d$-edges of $G$, $A^c$ is the set of 
$c$-edges of $G$. Then
$$G^A= G^{(A^c\Delta E(G))} = (G^{A^c})^*.$$
By Theorem~\ref{t.bipmed}, $G^{A^c}$ is bipartite, and, since the 
geometric dual of any bipartite graph is Eulerian (see Remark~\ref{r.btoe}), $(G^{A^c})^*$ is 
Eulerian as required.
\end{proof}

We note that the converse of Corollary~\ref{c.deu} is false and that determining exactly which subsets of edges of a plane graph give rise 
to Eulerian partial duals remains an open problem.

\begin{remark}\label{r.btoe}
The proof that the dual of a bipartite embedded graph is Eulerian is 
almost identical to the well-known proof of the special case for 
plane graphs: let $G$ be a bipartite graph, then every closed walk in 
$G$ is of even length (see \cite{ADH} for example). Therefore every 
closed walk about a face in any embedding of $G$ is of even length, 
and it follows that $G^*$ is Eulerian.
\end{remark}

\begin{remark}\label{r.knots}
The results on partial duals presented in this paper are intimately 
related to, and motivated by, the authors' work with N.~Virdee on the 
ribbon graphs of link diagrams in \cite{HMV}. In the context of knot 
theory, bipartite embedded graphs arise as the Seifert graphs of an 
oriented link diagram, and these embedded graphs are necessarily 
partial duals of plane graphs. This provides a link between the graph 
theory presented in Section~\ref{s.med} and knot theory. In fact, our 
knot theoretic results from \cite{HMV} on the characterization of 
Seifert graphs   suggested the formulation and proof of 
Theorem~\ref{t.bipmed} to us. Furthermore, the connection between the 
Tait graph and Seifert graph of a link diagram that was also studied 
in \cite{HMV}, led the authors to  conjecture that
$[(G\cup G^{*}) - (A^c  \cup  A^{*})]^{\du}$ and $G^A$
are isomorphic as abstract (but not necessarily embedded) graphs, 
which appears as Theorem~\ref{t.newpds} below. The results presented 
here therefore further illustrate the deep and fruitful connections 
between knot theory and graph theory.
\end{remark}

\section{The proof of Theorem~\ref{t.newpds}}\label{s.proof}

We will prove Theorem~\ref{t.newpds} by using the characterization of 
partially dual graphs in terms of a bijection between edge sets 
from~\cite{Mo4}. This extends the usual characterization of dual 
graphs in terms of maps between edge sets, which is due to Whitney 
for plane graphs \cite{Wh32}, and Edmonds for higher genus graphs 
\cite{Ed65}. We  need to introduce a little notation.

Suppose that $G$ and $H$ are graphs and $\varphi:E(G)\rightarrow 
E(H)$ is a bijection between their edge sets. Let $v\in V(G)$ and  $S\subseteq E(G)$.  Then we let  $S_v$ be the set of edges in $S$ which are incident with $v$. 
Then the set 
$\varphi(S_v)$ of edges in $H$ together with the vertices which are 
incident with these edges form a subgraph of $H$. This subgraph is denoted by $\varphi(S)_v$. 
Furthermore, we let $H_v$ denote the subgraph $\varphi(E(G))_v$.

\begin{definition}\label{d.ed}
Let $G$ and $H$ be graphs and $\varphi: E(G)\rightarrow E(H)$ be  a 
bijection. We say that $\varphi$ satisfies {\em Edmonds' Criteria} if 
\begin{enumerate}
\item edges $e,f\in E(G)$ belong to the same connected component if 
and only if $\varphi(e), \varphi(f)\in E(H)$ belong to the same 
connected component;
\item for each $v\in V(G)$, each component of $H_v$ is Eulerian;
\item for each $v\in V(H)$, each component of $G_v$ is Eulerian.
\end{enumerate}
\end{definition}

The significance of Edmonds' Criteria  is that they provide a 
characterization of geometric duality in terms of mappings between 
edge sets. Edmonds showed in  \cite{Ed65} that  $H\cong G^*$ for some 
cellular embedding of a graph $G$, if and only if there exists a 
bijection $\varphi:E(G)\rightarrow E(H)$ that satisfies Edmonds' 
Criteria and that $G$ and $H$  have the same number of isolated 
vertices. Edmonds' Theorem is an extension of Whitney's 
characterization of planar duals in terms of combinatorial duals (see 
\cite{Wh32}) to graphs embedded in an arbitrary surface.

In \cite{Mo4}, the second author extended Edmonds' and Whitney's 
Theorems to partial duals:
\begin{theorem}\label{t.main}
Two graphs $G$ and $H$ are partial duals if and only if there exists 
a bijection $\varphi: E(G) \rightarrow E(H)$, such that
\begin{enumerate}
\item \label{c1} $\left. \varphi\right|_{A}: A\rightarrow 
\varphi(A)$ satisfies Edmonds' Criteria for some subset $A\subseteq 
E(G)$.
\item \label{c2} If $v\in V(G)$ is incident with an edge in $A$, and 
if $e\in E(G)$ is incident with $v$, then $\varphi(e)$ is incident 
with a vertex of $\varphi(A)_v$. Moreover, if both ends of $e$ are incident 
with $v$, then both ends of $\varphi(e)$ are incident with vertices 
of $\varphi(A)_v$.
\item \label{c3} If $v\in V(G)$ is not incident with an edge in $A$, 
then there exists a vertex $v' \in V(H)$ with the property that $e\in 
E(G)$ is incident with $v$ if and only if $\varphi(e)\in E(H)$ is 
incident with $v'$. Moreover, both ends of $e$ are incident with $v$ 
if and only if both ends of $\varphi(e)$ are incident with $v'$.
\end{enumerate}
Furthermore, with $A$ as above $H\cong G^A $. 
\end{theorem}
We recall that $\varphi(A)_v $ is the subgraph of $H$ induced by the images of 
the edges from $A$ that are incident with $v$.

While the claim in the above theorem that  $H\cong G^A $ did not 
appear explicitly in \cite{Mo4}, it is an immediate consequence of  
the proof of Theorem~26  of \cite{Mo4}. We will use this 
characterization of partial duals to prove Theorem~\ref{t.newpds}.

\begin{figure}
    \subfigure[A vertex $v$ of $G$.]
    {
      \label{f.npd1}\raisebox{2mm}{\includegraphics[height=4cm]{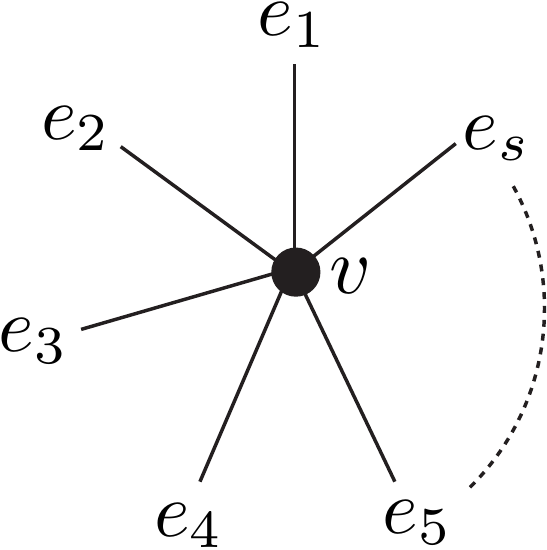} \hspace{5mm}}
    }
    \subfigure[A vertex $v$ of $G\cup G^{*}$.]
    {
      \label{f.npd2}\includegraphics[height=5cm]{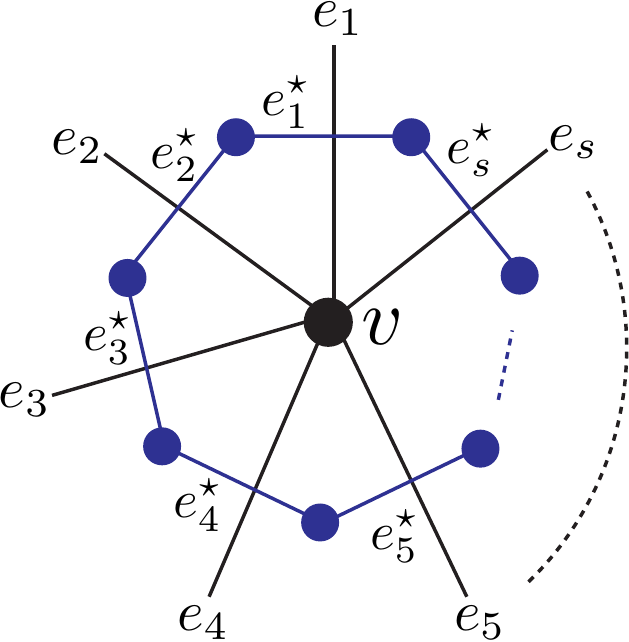}\hspace{5mm}
    } 
     \subfigure[Regions of $H$.]
     {
      \label{f.npd3}\includegraphics[height=5cm]{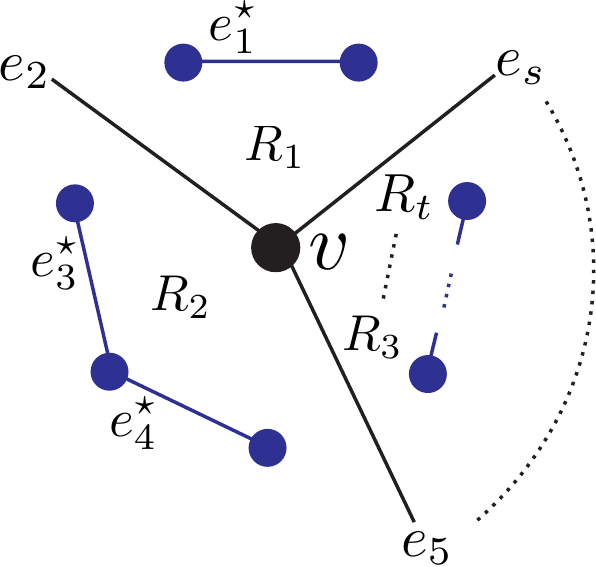}\hspace{5mm}
    }
    \subfigure[$H^{\du}$ and the mapping $\varphi$.]
    {
      \label{f.npd4}\raisebox{-5mm}{\includegraphics[height=6cm]{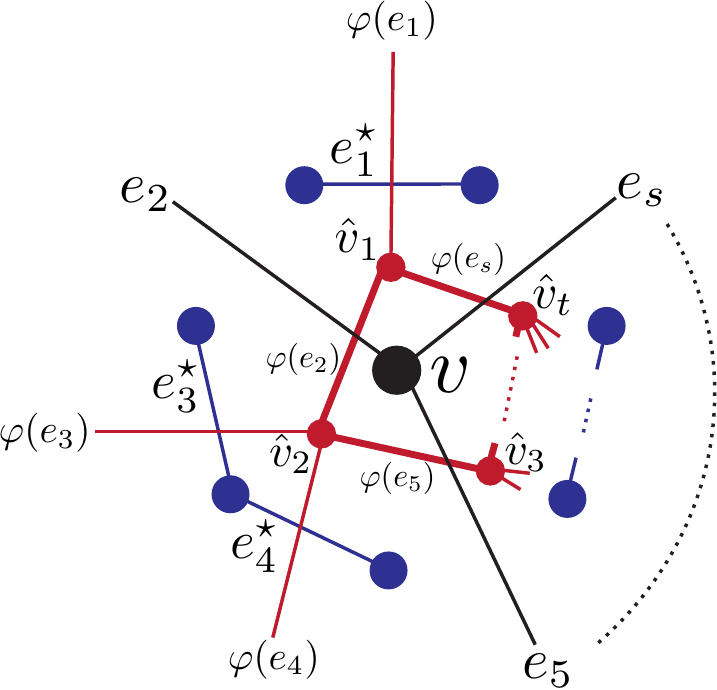}\hspace{5mm}}
    } 
  \caption{Figures used in the proof of Theorem~\ref{t.newpds}.}
  \label{f.npd}
\end{figure}

\begin{proof}[Proof of Theorem~\ref{t.newpds}]
To avoid clutter in the proof, we  set $H:= (G\cup G^{*}) - (A^c  
\cup  A^{*})$. We begin by defining a mapping $\varphi: E(G) 
\rightarrow E(H^{\du})$, and go on to show that this mapping 
satisfies the conditions of Theorem~\ref{t.main}.

Each edge $e\in E(G)$ is naturally identified with exactly one edge 
$e^{*}$ in $E(G^{*})$, and exactly one of $e$ or $e^{*}$ is in $H$. 
This gives a natural bijection $\alpha$ between $E(G)$ and $E(H)$. 
If we let $\beta$ denote the natural bijection between $E(H)$ and 
$E(H^{\du})$, then we obtain a natural bijection between $\varphi: 
E(G)\rightarrow E(H^{\du})$ by setting $\varphi:=\beta \circ 
\alpha$.

To prove the theorem it remains to show that $\varphi$ satisfies the 
conditions of Theorem~\ref{t.main} with the edge set $A$ given in the 
Theorem. For the first condition, we note that $A = E(G-A^c)$ and so, 
by the definition of $\alpha$, we have $\alpha(A)=E(G-A^c)$. Then, 
since $\beta$ is just the natural identification of edge sets of a 
graph and its dual, we have 
$$\varphi(A)=\beta(E(G-A^c))=E((G-A^c)^{\du}),$$
so
$$\left. \varphi\right|_{A}: E(G-A^c) \rightarrow E((G-A^c)^{\du}).$$
It is then readily verified that $\left. \varphi\right|_{A}$ 
satisfies Edmonds' Criteria. 

We will now show that $\varphi$ satisfies the remaining conditions of 
Theorem~\ref{t.main}. To do this we consider the construction of 
$H^{\du}$ from the embedded graph $G$ locally in the neighbourhood of 
a vertex $v$ of $G$. We begin with the immersed graph $G\cup G^{*}$. 
Let $v\in V(G)$, and let $e_1, e_2, \ldots ,e_s$ be the (not 
necessarily distinct) cyclically ordered edges of $G$ that are 
incident to $v$, where the cyclic order is chosen with respect to an 
arbitrary orientation of a neighbourhood of $v$. See 
Figure~\ref{f.npd1}. Let $e_i^{*}$ denote the unique edge in the 
subgraph $G^{*}$ of $G\cup G^{*}$ which intersects $e_i$, for each 
$i$. See Figure~\ref{f.npd2}. Note that the $e_i^{*}$s which arise 
need not be distinct. We will let $D_v$ denote the $s$-gon, together 
with the embedding of $v$ and its incident half-edges, which is 
obtained by cutting the immersed graph $G\cup G^{*}$ along the edges 
$e_1^{*}, \ldots ,e_s^{*}$ and their incident vertices. Next form $H$ 
by deleting all of the edges of the immersed graph $G\cup G^{*}$ 
which belong to $A^c\cup A^{*}$. The remaining edges of $G$ which are 
incident with $v$ divide $D_v$ into regions $R_1, \ldots , R_k$ in 
the following way: if no edges in $A$ are incident with $v$ then 
there is a single region $R_1$; if there are edges in $A$ that are 
incident with $v$, then cyclically order the regions $R_1, R_2, 
\ldots , R_t$ according to some orientation of $D_v$. See 
Figure~\ref{f.npd3}. Note that the regions $R_1, R_2, \ldots , R_t$ 
need not be distinct regions of the embedded graph $H$. Finally, form 
$H^{\du}$. There is a vertex $\hat{v}_k$ of $H^{\du}$ associated with 
each region $R_k$ (but the vertices $\hat{v}_1, \ldots , \hat{v}_t$ 
need not be distinct). See Figure~\ref{f.npd4}. In addition, notice 
that:
\begin{itemize}
\item if the edge $e_i$ of $H$ is adjacent to the region $R_k$, then 
$\varphi(e_i)$ is the edge which is incident with $\hat{v}_k$ and 
which intersects $e_i$ in the canonical immersion of $H\cup H^{\du}$; 
\item if the edge $e_i^{*}$ in $H$ is adjacent to the region $R_k$, 
then $\varphi(e_i)$ is the edge which is incident with $\hat{v}_k$ 
and which intersects $e_i^{*}$ in the canonical immersion of $H\cup 
H^{\du}$;
\item if $v$ is incident with an edge in $A$, then every $\hat{v}_k$ 
is in $\varphi(A)_v$.
\end{itemize}
We will use these observations (which are illustrated in 
Figure~\ref{f.npd4}) to verify that the map $\varphi : 
E(G)\rightarrow E(H^{\du})$ satisfies the remaining conditions of 
Theorem~\ref{t.main}. There are three  cases to consider.

\noindent {\bf Case 1:} Suppose that  $v$ is incident with an edge in 
$A$ and $e\in A$. In this case, by definition,  $\varphi(e)\in 
\varphi(A)_v$ and so both ends of $\varphi(e)$ are in 
$\varphi(A)_v$.

\noindent {\bf Case 2:} Suppose that  $v$ is incident with an edge in 
$A$ and $e\notin A$. In this case, we have that $\varphi(e)=e_i^{*}$, 
for some $i$. Then $e_i^{*}$ is adjacent to a region $R_k$, and hence 
$\varphi(e)$ is incident with $\hat{v}_k$, and therefore incident 
with a vertex of $\varphi(A)_v$.

If $e$ is a loop then $\varphi(e)=e_i^{*} = e_j^{*}$, for some $i$ 
and $j$. Then $e_i^{*}$ is adjacent to a region $R_k$, and $e_i^{*}$ 
is adjacent to a region $R_l$. This means that $\varphi(e)$ is 
incident with $\hat{v}_k$ and with $\hat{v}_l$ (with incidence being 
counted twice if $\hat{v}_k = \hat{v}_l$). Thus, both ends of 
$\varphi(e)$ are incident with $\varphi(A)_v$.

\noindent {\bf Case 3:} Suppose that $v$ is not incident with an edge 
in $A$. In this case there is only one region $R_1$, and the vertex 
$\hat{v}_1$ is the vertex $v'$ required by Theorem~\ref{t.main}. To 
see why this is the case, suppose that $e\in E(G)$ is incident with 
$v$. Then $\varphi(e)=e_i^{*}$, for some $i$, and, as before, 
$e_i^{*}$ is adjacent to a region $R_1$, so $\varphi(e)$ is incident 
with $\hat{v}_1$, as required. If, in addition, $e$ is a loop, then 
$\varphi(e)=e_i^{*} = e_j^{*}$, for some $i$ and $j$. Then since 
$e_i^{*}$ and $e_j^{*}$ are both adjacent to $R_1$ both ends of 
$\varphi(e)$ are incident to  $\hat{v}_1$.

Thus we have shown that $\varphi$ satisfies the conditions of 
Theorem~\ref{t.main}, and so $H\cong G^A$ as required.

\end{proof}

\end{document}